\documentclass{amsart}

\usepackage{latexsym,amscd,amsmath,verbatim,amssymb}

\theoremstyle{plain}
\newtheorem{theorem}{Theorem}[section]
\newtheorem{corollary}[theorem]{Corollary}

\newtheorem{lemma}[theorem]{Lemma}
\newtheorem*{theorem 1.1}{Theorem 1.1}
\newtheorem*{theorem 2.2}{Theorem 2.2}
\newtheorem*{theorem 3.3}{Theorem 3.3}
\newtheorem*{corollary 3.4}{Corollary 3.4}

\theoremstyle{remark}
\newtheorem{remark}[theorem]{Remark}
\newtheorem{notation}[theorem]{Notation}
\newtheorem*{remark 0.3}{Observations on independence 0.3}

\theoremstyle{definition}
\newtheorem{definition}[theorem]{Definition}
\newtheorem{example}[theorem]{Example}

\numberwithin{equation}{section}

\newcommand{\ann}{\operatorname{Ann}}
\newcommand{\ass}{\operatorname{Ass}}
\newcommand{\spec}{\operatorname{Spec}}
\newcommand{\ar}{\operatorname{AR}}

\begin{document}

\title[Primary Decomposition]{Primary Decomposition: Compatibility,
  Independence and Linear Growth}
\author{Yongwei Yao}
\address{Department of Mathematics, University of Kansas, Lawrence,
  Kansas 66045}   
\email{yyao@math.ukans.edu}
\begin{abstract}
For finitely generated modules $N \subsetneq M$ over a Noetherian ring 
$R$, we study the following properties about primary
decomposition: (1) The Compatibility property, which says that if 
$\ass (M/N)=\{ P_1, P_2, \dots , P_s\}$ and $Q_i$ is a $P_i$-primary
component of $N \subsetneq M$ for each $i=1,2,\dots,s$, then $N =Q_1
\cap Q_2  \cap \cdots \cap Q_s$; (2) For a given subset $X=\{ P_1,
P_2,  \dots , P_r \} \subseteq \ass(M/N)$, $X$ is an open subset of
$\ass(M/N)$ if and only if the intersections $Q_1 \cap
Q_2\cap \cdots \cap Q_r= Q_1' \cap Q_2' \cap \cdots \cap Q_r'$ for all
possible $P_i$-primary components $Q_i$ and $Q_i'$ of $N\subsetneq M$;
(3) A new proof of 
the `Linear Growth' property, which says that for any fixed ideals
$I_1, I_2, \dots, I_t$ of $R$, there exists a $k \in \mathbb N$ such that
for any $n_1, n_2, \dots, n_t \in \mathbb N$ there exists a primary
decomposition of $I_1^{n_1}I_2^{n_2}\cdots I_t^{n_t}M \subset M$ such
that every $P$-primary component $Q$ of that primary decomposition
contains $P^{k(n_1+n_2+\cdots+n_t)}M$. 
\end{abstract}
\maketitle

\section*{0. Introduction} 

Throughout this paper $R$ is a Noetherian ring and $M \neq 0$
is a finitely generated $R$-module unless stated otherwise
explicitly. Let $N 
\subsetneq  M$ be a proper $R$-submodule of $M$.  By primary
decomposition $N= Q_1 \cap Q_2  \cap \cdots \cap Q_s$ of $N$ in $M$,
we  always mean an irredundant and minimal primary decomposition,
where $Q_i$ 
is a $P_i$-primary submodule of $M$, i.e. $\ass(M/Q_i) = \{P_i\}$,
for  each $i=1, 2, \dots , s$,
unless mentioned otherwise explicitly. Then $\ass (M/N)=\{ P_1, P_2,
\dots , P_s\}$ and we say that $Q_i$ is a \mbox{$P_i$-primary}
component  of $N$ in $M$. As a subset of $\spec (R)$ with
the {\it Zariski} topology, $\ass (M/N)$ inherits a topology
structure. For an ideal $I$ in $R$, we  use  $(N :_M I^{\infty})$ to
denote $\cup_i(N :_M I^i)$.

\begin{notation}Let $N \subsetneq  M$ be finitely
generated $R$-modules. For every $P \in \ass(M/N)$, we use
$\Lambda_P(N \subsetneq M)$, or $\Lambda_P$ if the
$R$-modules $N \subsetneq  M$ are clear from the context,  
to denote the set of all possible
$P$-primary components of $N$ in $M$. 
\end{notation}

We know that if $P \in \ass(M/N)$ is an embedded prime ideal, then
$\Lambda_P(N\subsetneq M)$ contains more than one element. (Also see the
passage following Theorem 2.2 and the reference to \cite{HRS}.)
Suppose that 
$N= Q_1 \cap Q_2  \cap \cdots \cap Q_s$ is a primary decomposition of
$N\subsetneq M$ such that $Q_i \in \Lambda_{P_i}$ for $i = 1,2,\dots, s$.
Then if we choose a $P_i$-primary submodule $Q_i'$ of $M$ such that
$N\subseteq Q_i' \subseteq Q_i$ for each $i=1,2,\dots,s$, we get a
primary decomposition $N= Q_1' \cap Q_2'  \cap \cdots \cap Q_s'$ of $N
\subsetneq M$. For example we may choose $Q_i'=\ker (M \to
(M/(P_i^{n_i}M+N))_{P_i})$ for all $n_i\gg 0$ to get primary
decompositions $N=\cap_{1\le i\le s}\ker(M \to
(M/(P_i^{n_i}M+N))_{P_i})$ for all $n_i\gg 0$. But given an arbitrary
$Q_i'' \in 
\Lambda_{P_i}$ for each $i=1,2,\dots,s$, we do not know \it a
priori \rm if $N=Q_1'' \cap Q_2''  \cap \cdots \cap Q_s''$. This
compatibility question is answered positively in Theorem 1.1:  

\begin{theorem 1.1}[Compatibility] Let $N\subsetneq M$ be finitely
generated $R$-modules  and
$\ass (M/N)=\{ P_1, P_2, \dots , P_s\}$. Suppose that for each $i=1,2,
\dots ,s$, $Q_i$ is a \mbox{$P_i$-primary} component of $N$ in $M$, i.e. $Q_i
\in \Lambda_{P_i}$. Then $N=Q_1 \cap Q_2\cap \cdots \cap Q_s$,
which is necessarily an irredundant and minimal primary decomposition.
\end{theorem 1.1}

\begin{definition} Let $N \subsetneq M$ be  finitely
generated  $R$-modules and $X$ a subset of $\ass(M/N)$, say 
$ X =\{ P_1,  P_2,  \dots , P_r \} \subseteq \ass (M/N) = \{ P_1,
\dots ,  P_r, P_{r+1}, \dots, P_s\}$. We say that the primary 
decompositions of $N$ in $M$ are independent over $X$, or
$X$-independent,
if for any two primary decompositions, say, $N= Q_1 \cap Q_2
\cap \cdots \cap Q_s = Q'_1 \cap Q'_2  \cap \cdots \cap Q'_s$, of
$N\subset M$ such that $\{Q_i,Q'_i\} \subseteq \Lambda_{P_i}(N \subset
M)$ for $i=1,2,\dots,s$, we have $Q_1 \cap Q_2\cap \cdots \cap
Q_r=Q'_1 \cap Q'_2 \cap \cdots \cap Q'_r$. In this case, we denote the
invariant intersection by $Q_X(N \subset M)$, or $Q_X$ if $N \subset
M$ is clear from the context.
\end{definition}

It is well-known that primary decompositions are independent over open
subsets of $\ass(M/N)$. (See Observations 0.3 below.) Actually it
turns out that independence property characterizes open subsets of
$\ass(M/N)$:

\begin{theorem 2.2} Let $N \subsetneq M$ be finitely generated $R$-modules
and $X \subseteq \ass(M/N)$ be a subset of $\ass(M/N)$. Then the primary
decompositions of $N$ in $M$ are independent over $X$ if
and only if $X$ is an open subset of $\ass(M/N)$.
\end{theorem 2.2}

In Section 3 we use \it Artin-Rees numbers \rm to prove the following:

\begin{theorem 3.3} Let $R$ be a Noetherian ring,  $M$ a finitely
generated $R$-module and $I_1, I_2, \dots, I_t$ ideals  of
$R$. Then there exists a $k \in \mathbb N$ such that for all $n_1, n_2,
\dots , n_t \in \mathbb N$ and 
for all ideals $J \subset R$,  $(J^{k|\underline
n|}M+I_1^{n_1}I_2^{n_2}\cdots I_t^{n_t}M) \cap 
(I_1^{n_1}I_2^{n_2}\cdots I_t^{n_t}M :_M  J^{\infty})=
I_1^{n_1}I_2^{n_2}\cdots I_t^{n_t}M$,
where $|\underline n|:=n_1+n_2+\cdots +n_t$. 
\end{theorem 3.3} 

As a corollary of Theorem 3.3, we have a new proof of the `Linear
Growth' property, which was first proved by I. Swanson \cite{Sw} and
then  by R. Y. Sharp using different methods and in a more general
situation \cite{Sh2}:

\begin{corollary 3.4}[Linear Growth; \cite{Sw} and \cite{Sh2}]
Let $R$ be a Noetherian ring,   $M$  a
finitely generated $R$-module and $I_1, I_2, \dots, I_t$ ideals  of
$R$. Then there exists a $k \in \mathbb N$ such that for any $n_1, n_2,
\dots , n_t \in \mathbb N$, there exists a primary decomposition of $
I_1^{n_1}I_2^{n_2}\cdots I_t^{n_t}M \subseteq M$  
\[ 
I_1^{n_1}I_2^{n_2}\cdots I_t^{n_t}M = Q_{\underline n_1} \cap
Q_{\underline n_2} \cap \cdots \cap  Q_{\underline n_{r_{\underline
n}}},
\]
where the $Q_{\underline n_i}$'s are $P_{\underline n_i}$-primary
components of the primary decomposition 
such that $P_{\underline n_i}^{k|\underline n|}M \subseteq
Q_{\underline n_i}$ for all $i=1,2,\dots,  
{r_{\underline n}}$, where $\underline n=(n_1, n_2,
\cdots, n_t)$ and $|\underline n|=n_1+n_2+\cdots +n_t$.  
\end{corollary 3.4}   

Before ending this introduction section, we make the following
well-known observations, which is to the effect of saying that primary
decompositions are independent over open subsets. 

\begin{remark 0.3} Suppose $N= Q_1  \cap Q_2 \cap
\cdots \cap Q_s$ is a primary decomposition of $N$ in a finitely
generated $R$-module $M$ such that $Q_i$ is $P_i$-primary for each
$i=1,2,\dots,s$. 
\begin{enumerate}  
\item For any ideal $I \subseteq R$, the intersection
$ \cap_{I \not \subseteq P_i} Q_i =  (N:_M I^{\infty})$ is 
independent of the particular primary decomposition of $N$ in
$M$. (cf. D. Eisenbud \cite{Ei}, page 101, Proposition 3.13.) This
means that the primary decompositions of $N \subsetneq M$ are independent
over $X = \{ P \in \ass(M/N) \,|\, I \not \subseteq P\}$ and
$Q_X=(N:_M I^{\infty})$.  
\item Alternatively, for any multiplicatively  
closed set $W\subset R$, the intersection
$ \cap_{P_i \cap W = \emptyset} Q_i = \ker (M \to (M/N)_W)$
is independent of the particular primary decomposition.
(cf. D. Eisenbud \cite{Ei}, page 113, Exercise  3.12.) That is to say
that the primary decompositions of $N \subsetneq M$ are independent
over $Y=\{P \in \ass (M/N) \,|\, P \cap W = \emptyset \}$ and
$Q_Y=\ker (M \to (M/N)_W)$.
\end{enumerate}
\end{remark 0.3}

\section{Compatibility}

The main theorem in this section is to show that all the primary
components of $R$-modules $N \subsetneq M$ are totally compatible in
forming the primary decompositions of $N \subsetneq M$.

\begin{theorem 1.1}[Compatibility] Let $N\subsetneq M$ be finitely
generated $R$-modules  and
$\ass (M/N)=\{ P_1, P_2, \dots , P_s\}$. Suppose that for each $i=1,2, \dots
,s$, $Q_i$ is a \mbox{$P_i$-primary} component of $N$ in $M$, i.e. $Q_i \in
\Lambda_{P_i}(N \subsetneq M)$. Then $N=Q_1 \cap Q_2\cap \cdots \cap Q_s$,
which is necessarily an irredundant and minimal primary decomposition.
\end{theorem 1.1}

\begin{proof} We induct on $s$, the cardinality of $\ass(M/N)$. 

If $s=1$, then $N=Q_1$ and the claim is trivially true. 

Suppose $s \ge 2$. By rearranging the order of $P_1, P_2, \dots , P_s$,
we may assume that $P_s$ is a maximal prime ideal in $\ass(M/N)$.
Since $Q_i \in \Lambda_{P_i}$ for $i=1,2,\dots, s$, we can find $s$ 
specific primary decompositions
\[ 
N=Q_{(i,1)}\cap Q_{(i,2)}\cap \cdots \cap Q_{(i,i)}\cap \cdots \cap
Q_{(i,s)}, \qquad \text {for } i=1,2,\dots,s,
\]
where $Q_{(i,j)} \in \Lambda_{P_j}$ and $Q_{(i,i)}=Q_i$ for all
$i,j=1,2,\dots, s$. Let $W=R \setminus \cup_{1 \le i \le s-1}P_i$.
By Observation 0.3(2) and our assumption on $P_s$,
we know that the primary decompositions of $N\subsetneq M$ is independent
over $X=\{ P \in \ass(M/N) \,|\, P \cap W= \emptyset \}=\{ P_1, P_2,
\dots , P_{s-1}\}$ with $Q_X=\ker(M \to (M/N)_W)$. That is to say that 
\[
Q_X=\ker(M \to (M/N)_W)= Q_{(i,1)}\cap Q_{(i,2)}\cap \cdots \cap
Q_{(i,s-1)}, \qquad \text {for } i=1,2,\dots,s,
\]
are all primary decompositions of $Q_X \subsetneq M$ and in particular
$Q_i=Q_{(i,i)}\in \Lambda_{P_i}(Q_X \subset M)$ for $i=1,2,
\dots,s-1$. Since the
cardinality of $\ass(M/Q_X)$ is $s-1$, we use the induction hypothesis
to see that 
\[
Q_X=Q_1 \cap Q_2\cap \cdots \cap Q_{s-1}.
\]
But we already know that $Q_X= Q_{(s,1)}\cap Q_{(s,2)}\cap \cdots \cap
Q_{(s,s-1)}$ by the $X$-independence of primary decompositions of
$N\subsetneq M$. Hence we have
\begin{eqnarray*}
N&=&Q_{(s,1)}\cap Q_{(s,2)}\cap \cdots \cap Q_{(s,s-1)}\cap Q_{(s,s)}\\
 &=&Q_X \cap Q_s\\
 &=&Q_1 \cap Q_2\cap \cdots \cap Q_{s-1}\cap Q_s. \end{eqnarray*}
\end{proof}

\addtocounter{theorem}{1}

\begin{remark} In \cite[Chapter IV]{Bo}, the
notion of primary decomposition is generalized to not necessarily
finitely generated modules over not necessarily Noetherian rings. Let
$R$ be a (not necessarily Noetherian) ring and $M$ be a (not
necessarily finitely generated) $R$-module. A prime ideal $P$ of $R$
is said to be \it weakly associated \rm with $M$ if there exists an $x
\in M$ such that $P$ is minimal over the ideal $\ann(x)$ and we denote
by $\ass_f(M)$ the set of prime ideals weakly associated with $M$ 
(cf. \cite[page 289, Chapter IV, $\S$ 1, Exercise 17]{Bo}.)  
We say that an element $r
\in R$ is nearly nilpotent on $M$ if for any $x \in M$, there exists
an $n(x) \in \mathbb N$, such that $r^{n(x)}x=0$ (cf. \cite[{page 267,
Chapter IV, $\S$ 1.4, Definition 2}]{Bo}.) Then for any 
$R$-submodule $N$ of $M$, we define
$r_M(N):=\{r \in R \,|\, r \text { is nearly nilpotent on } M/N
\}$ (cf. \cite[page 292, Chapter IV, $\S$ 2, Exercise 11]{Bo}.) A
$R$-submodule $Q$ of $M$ is said to be $P$-primary in $M$ if
$\ass_f(M/Q) = \{ P \}$, which is equivalent to the statement that
every $r \in R$ is either a non-zerodivisor or nearly nilpotent on
$M/Q$, and in this case we have $r_M(Q)=P)$ (cf. \cite[page 292,
Chapter IV, $\S$ 2, Exercise 12(a)]{Bo}.) 
Then we say that a $R$-submodule $N$ has a primary decomposition in
$M$ if there exist $P_i$-primary submodules $Q_i \subset M$,
$i=1,2,\dots, s$, such that $N=Q_1 \cap Q_2\cap \cdots \cap
Q_s$ (cf. \cite[page 294, Chapter IV, $\S$ 2, Exercise 20]{Bo}.)
Again we always assume primary decompositions to be 
irredundant and minimal (i.e. reduced) if they exist. If $N$ has
primary decompositions in $M$, then Observation 0.3(2) still holds 
(replace $\ass(M/N)$ by $\ass_f(M/N)$.) Therefore the proof of
compatibility, i.e. Theorem 1.1, also applies to the the case
where $N\subsetneq M$ are not necessarily finitely generated
$R$-modules over a not necessarily Noetherian ring $R$ as long as the
primary decompositions exist.
\end{remark}

\section{Independence over open subsets of $\ass(M/N)$} 

Because of the compatibility property, i.e. Theorem 1.1, we have an
equivalent statement to the definition of $X$-independence.

\begin{lemma} Let  $N \subsetneq M$ be  finitely
generated  $R$-modules and
$ X =\{ P_1,  P_2,  \dots , P_r \} \linebreak \subseteq \ass (M/N) =
\{ P_1, 
P_2, \dots ,  P_r, P_{r+1}, \dots, P_s\}$. Then the following are
equivalent: 
\begin{enumerate}
\item The primary decompositions of $N$ in $M$ are
independent over $X$;
\item For any $Q_i$ and $Q'_i$ in $\Lambda_{P_i}$, where
$i=1,2,\dots, r$, the equality \linebreak
 $Q_1 \cap Q_2\cap \cdots \cap
Q_r=Q'_1 \cap Q'_2 \cap \cdots \cap Q'_r$ holds.
\end{enumerate}
\end{lemma}

It turns out that the independence observed in
Observations 0.3 actually exhausts all the possibilities. 

\begin{theorem} Let $N \subsetneq M$ be finitely generated $R$-modules
and $X \subseteq \ass(M/N)$ be a subset of $\ass(M/N)$. Then the primary
decompositions of $N$ in $M$ are independent over $X$ if
and only if $X$ is an open subset of $\ass(M/N)$.
\end{theorem}

\begin{proof} Without loss of generality we assume $N=0$.

The ``if'' part is just Observation 0.3(1). To prove the
``only if'' part, it suffices to show $X$ is stable under
specialization since $\ass(M/N)=\ass(M)$ is finite. Let $P$ be an arbitrary
prime ideal in $X \subseteq \ass(M/N)$. 
All we need to show is that for any $P'\in \ass(M)$ such that $P'
\subset P$, we have $P' \in X$.

Say $X= \{P= P_1,  P_2,  \dots , P_t, P_{t+1},\dots , P_r
\}$  such that $P_i \subseteq P$ for $i=1,2,\dots,t$ and $P_i \not
\subseteq P$  for $i=t+1,\dots,r$. Let $X_P:=X \cap \ass(M_P) = \{P_P=
(P_1)_P, (P_2)_P, \linebreak[2] \dots , (P_t)_P \}$. We first show
that the primary 
decompositions of $0 \subsetneq M_P$ are independent over $X_P$:
For any $L_i \in \Lambda_{(P_i)_P}(0 \subsetneq 
M_P), i=1,2,\dots,t,$ let $Q_i$ be the the full pre-image of $L_i$
under the map $M \to M_P$. Then choose $Q_i \in \Lambda_{P_i}(0\subsetneq
M)$ for $i= t+1,\dots,r$. Then it is easy to see that $(Q_1\cap Q_2
\cap \cdots \cap Q_r)_P= L_1\cap L_2 \cap \cdots \cap
L_t$. Then the $X$-independence assumption implies that the primary
decompositions  of $0 \subsetneq M_P$ are independent over $X_P = X \cap
\ass(M_P)$.

Hence by replacing $M$ with $M_P$ we may assume that $(R, P)$ is local
with the maximal ideal $P$ and 
$P \in X= \{P= P_1,  P_2,  \dots , P_t\} \subseteq \ass(M)$. In this
case to prove that $X$ is stable under specialization is simply to
prove that $X = \ass (M)$. 
For each $i=1,2,\dots,t$, choose a $P_i$-primary component $Q_i$ of
$0 \subsetneq M$. There exists a $k \in \mathbb N$ such that $P^kM
\subseteq Q_1$ and therefore $P^nM \in \Lambda_P$ for all $n\ge k$. 
Set $L = Q_2 \cap Q_3 \cap \cdots \cap Q_t$.
Then by Lemma 2.1 the assumption that the primary decompositions of
$0$ in $M$ are independent over $X$  simply
means that $Q_1 \cap L = P^nM \cap L$ for all $n\ge k$, which implies
$Q_1 \cap L = 0$ by Krull Intersection Theorem. This forces $0 =Q_1 
\cap Q_2 \cap \cdots \cap Q_t$ to be a primary decomposition of $0$ in
$M$. In particular it means that $\ass(M) = \{P= P_1,  P_2,  \dots,
P_t \} = X$. 
\end{proof}

In particular, if $P \in \ass(M/N)$ is not minimal over $\ann(M/N)$,
then the \mbox{$P$-primary} components of $N$ in $M$ are not
unique. In fact, in \cite{HRS},  W. Heinzer, 
L. J. Ratliff, Jr. and K. Shah showed that if $P \in \ass(M/N)$ is an
embedded prime ideal, then there are infinitely many maximal
$P$-primary components of $N$ in $M$ with respect to containment. See
\cite{HRS} and their following papers for more information about the
embedded primary components.

\section{`Linear Growth' property}

In this section we give a new proof of `Linear Growth' property using
\it \mbox{Artin-Rees} numbers \rm and compatibility. `Linear Growth'
property was first proved by I. Swanson \cite{Sw} and then 
by R. Y. Sharp using different methods and in a more general situation
\cite{Sh2}.  

We first give a definition of \it Artin-Rees numbers, \rm 
$\ar(J,N \subseteq M)$, of a pair of finitely generated $R$-modules $N
\subseteq M$ with respect to an ideal $J$ of $R$. These numbers have been
studied in \cite{Hu}, where a set of ideals is considered instead
of one single ideal.

\begin{definition} Let $N \subseteq M$ be finitely generated
$R$-modules over a Noetherian ring $R$ and $J$ an ideal of $R$. We
define $\ar(J,N \subseteq M):=\min \{\,k \,|\, J^nM \cap N \subseteq
J^{n-k}N \text{ for all } n \ge k \,\}$.
\end{definition}

\begin{remark} If $K \subseteq L \subseteq M$, then $\ar(J,K\subseteq
M) \le \ar(J,K\subseteq L)+\ar(J,L\subseteq M)$. If $J^nM \subseteq
N$ for some $n$, then $\ar(J, N \subseteq M) \le n$.
\end{remark}

\begin{theorem} Let $R$ be a Noetherian ring,  $M$ a finitely
generated $R$-module and $I_1, I_2, \dots, I_t$ ideals  of
$R$. Then there exists a $k \in \mathbb N$ such that for all $n_1, n_2,
\dots , n_t \in \mathbb N$ and 
for all ideals $J \subset R$, 
\begin{eqnarray*}
I_1^{n_1}I_2^{n_2}\cdots I_t^{n_t}M &\supseteq&
J^{k|\underline n|}M \cap (I_1^{n_1}I_2^{n_2}\cdots 
I_t^{n_t}M :_M  J^{\infty}), \qquad \text {i.e.} \\
I_1^{n_1}I_2^{n_2}\cdots I_t^{n_t}M &=&
(J^{k|\underline n|}M+I_1^{n_1}I_2^{n_2}\cdots I_t^{n_t}M) \cap
(I_1^{n_1}I_2^{n_2}\cdots I_t^{n_t}M :_M  J^{\infty}), \end{eqnarray*}
where $|\underline n|:=n_1+n_2+\cdots +n_t$. 
\end{theorem}

\begin{proof}It is enough to prove the Theorem for
\begin{eqnarray*}
\mathcal{R}&=&R[I_1T_1,I_2T_2,\dotsc, I_tT_t,
T_1^{-1},T_2^{-1},\dots,T_t^{-1}], \\
\mathcal{M}&=& \bigoplus_{n_1, n_2, \dotsc, n_t \in \mathbb Z}
I_1^{n_1}I_2^{n_2}\cdots  I_t^{n_t}MT_1^{n_1}T_2^{n_2}\cdots
T_t^{n_t},\\ 
\mathcal{I}_i&=&T_i^{-1} \mathcal{R} \qquad \text{for each }
i=1,2,\dots,t, \qquad \text{and}\\ 
\mathcal{J}&=&J\mathcal{R}. 
\end{eqnarray*}
That is because if we contract the result for $\mathcal{R}$ back to $R$, we
get the desired result.
Hence without loss of generality we assume $I_i=(x_i)$ is generated by
a $M$-regular element $x_i\in R$ for each
$i=1,2,\dots,t$.
The same technique is also used in \cite{Sw} and \cite{Sh2}.

And it also suffices to prove the Theorem for one
fixed ideal $J$. The reason is for every $J$ in $R$, we have 
\[
J \subseteq J':= \bigcap_{\substack{P \in Y\\J \subseteq P}} P,
\quad
\text{where}\quad Y= \bigcup_{(n_1, n_2, \dots , n_t) \in \mathbb Z^t}
\ass(M/I_1^{n_1}I_2^{n_2}\cdots I_t^{n_t}M)
\] 
and, furthermore, there are only finitely many such $J'$ to deal
with since the set
$Y=\cup_{(n_1,n_2,\cdots,n_t)}\ass(M/I_1^{n_1}I_2^{n_2}\cdots 
I_t^{n_t}M)$ is finite. (cf. \cite[page 125, Lemma 13.1]{Mc}) 

For each $i=1,2,\dots,t$, let $N_i =x_iM:_M J^{\infty} \subseteq
M$, $k_i'=\ar(J,N_i \subseteq M)$ and $k_i''$ be such that
$J^{k_i''}N_i \subseteq x_iM$. Then $\ar (J,x_iM \subseteq N_i) \le
k_i''$. 

Let $k'=\max
\{ k_i'  \,|\, 1\le i \le t\,\}$, $k''=\max \{ k_i''  \,|\, 1\le i \le
t\,\}$ and $k=k' + k''$.
It is easy to see by the Remark 3.2
$\ar(J,x_iM \subseteq M) \le k_i'+k_i'' \le k$ for all
$i=1,\dots,t$. Since each 
$x_i$ is regular on $M$, we have $\ar (J, x_1^{m_1}x_2^{m_2}\cdots
x_{i-1}^{m_{i-1}}x_i^{m_i+1}x_{i+1}^{m_{i+1}}\cdots x_t^{m_t}M \subseteq 
x_1^{m_1}x_2^{m_2} \cdots x_t^{m_t}M)= \ar (J,x_iM \subseteq M) \le k$
because of the $R$-linear isomorphism $M \cong x_1^{m_1}x_2^{m_2}
\cdots x_t^{m_t}M$ induced by multiplication by 
$x_1^{m_1}x_2^{m_2} \cdots x_t^{m_t}$. Therefore 
we have $\ar(J,x_1^{n_1}x_2^{n_2}\cdots x_t^{n_t}M \subseteq M) 
\le k(n_1+n_2+\cdots +n_t) =k|\underline n|$ by the same Remark
3.2 applied to the filtration 
\[
x_1^{n_1}x_2^{n_2}\cdots x_t^{n_t}M \subseteq
x_1^{n_1-1}x_2^{n_2}\cdots x_t^{n_t}M \subseteq 
\cdots \subseteq 
x_t^2M \subseteq x_tM \subseteq M
\]
of $x_1^{n_1}x_2^{n_2}\cdots x_t^{n_t}M \subseteq
M$ so that each quotient 
is isomorphic to $M/x_iM$ for some $i=1,2,\dots,t$. 

We prove the Theorem by induction on $|\underline n|=n_1+n_2+\cdots
+n_t$. If $|\underline n|=0$, the claim is trivially true.
 
Now suppose $|\underline n| \ge 1$. By symmetry we assume $n_1 \ge 1$. 
Notice, by induction hypothesis,
\begin{equation*}
\begin{split}
J^{k|\underline n|}M \cap &(x_1^{n_1}x_2^{n_2}\cdots x_t^{n_t}M:_M
J^{\infty})\\  
\subseteq &J^{k(|\underline n|-1)}M \cap (x_1^{n_1-1}x_2^{n_2}\cdots
x_t^{n_t}M:_M J^{\infty})\\ 
\subseteq &x_1^{n_1-1}x_2^{n_2}\cdots
x_t^{n_t}M. 
\end{split}\tag{$*$}
\end{equation*} 
Therefore, using the definition of integers  $k,\ k',\ k''$  and  the
fact that   
\begin{equation*}
\begin{split}
\ar\Big(J,\ (x_1^{n_1}x_2^{n_2}\cdots x_t^{n_t}M:_{x_1^{n_1-1}x_2^{n_2}\cdots
x_t^{n_t}M} &J^{\infty}) \subseteq
x_1^{n_1-1}x_2^{n_2}\cdots x_t^{n_t}M\Big)\\&=\ar(J,\ x_1M:_M J^{\infty}
\subseteq M) \qquad \text{and}\\  
(x_1^{n_1}x_2^{n_2}\cdots x_t^{n_t}M:_{x_1^{n_1-1}x_2^{n_2}\cdots
x_t^{n_t}M} J^{\infty})/&x_1^{n_1}x_2^{n_2}\cdots x_t^{n_t}M\\ &\cong
(x_1M :_M J^{\infty})/x_1M, 
\end{split}
\end{equation*}
we have,  
\begin{equation*}
\begin{split}
J^{k|\underline n|}M \cap &(x_1^{n_1}x_2^{n_2}\cdots x_t^{n_t}M:_M
J^{\infty})\\ 
=&(x_1^{n_1-1}x_2^{n_2}\cdots x_t^{n_t}M) \cap J^{k|\underline n|}M \cap
(x_1^{n_1}x_2^{n_2} \cdots
x_t^{n_t}M:_M  J^{\infty}) \qquad \text{by $(*)$} \\
= &(x_1^{n_1-1}x_2^{n_2}\cdots x_t^{n_t}M) \cap J^{k|\underline n|}M
\cap (x_1^{n_1}x_2^{n_2}\cdots x_t^{n_t}M:_{x_1^{n_1-1}x_2^{n_2}\cdots
x_t^{n_t }M} J^{\infty}) \\
= &\big(J^{k|\underline n|}M \cap (x_1^{n_1-1}x_2^{n_2} \cdots
x_t^{n_t}M)\big) 
\cap (x_1^{n_1}x_2^{n_2}\cdots x_t^{n_t}M:_{x_1^{n_1-1}x_2^{n_2}\cdots
x_t^{n_t }M} J^{\infty}) \\
\subseteq &\big(J^k(x_1^{n_1-1}x_2^{n_2}\cdots x_t^{n_t}M)\big)
\cap (x_1^{n_1}x_2^{n_2}\cdots x_t^{n_t}M:_{x_1^{n_1-1}x_2^{n_2}\cdots
x_t^{n_t}M} J^{\infty}) \\ 
\subseteq &J^{k''}(x_1^{n_1}x_2^{n_2}\cdots
x_t^{n_t}M:_{x_1^{n_1-1}x_2^{n_2}\cdots
x_t^{n_t}M} J^{\infty})\\
\subseteq &x_1^{n_1}x_2^{n_2}\cdots x_t^{n_t}M. 
\end{split}\end{equation*}
\end{proof}

\begin{corollary 3.4} [Linear Growth; \cite{Sw} and \cite{Sh2}]
Let $R$ be a Noetherian ring,   $M$  a
finitely generated $R$-module and $I_1, I_2, \dots, I_t$ ideals  of
$R$. Then there exists a $k \in \mathbb N$ such that for any $n_1, n_2,
\dots , n_t \in \mathbb N$, there exists a primary decomposition of
\linebreak[4] $I_1^{n_1}I_2^{n_2}\cdots I_t^{n_t}M \subseteq M$  
\[ 
I_1^{n_1}I_2^{n_2}\cdots I_t^{n_t}M = Q_{\underline n_1} \cap
Q_{\underline n_2} \cap \cdots \cap  Q_{\underline n_{r_{\underline
n}}},
\]
where the $Q_{\underline n_i}$'s are $P_{\underline n_i}$-primary
components of the primary decomposition 
such that $P_{\underline n_i}^{k|\underline n|}M \subseteq
Q_{\underline n_i}$ for all $i=1,2,\dots, 
{r_{\underline n}}$, where $\underline n=(n_1, n_2,
\cdots, n_t)$ and $|\underline n|=n_1+n_2+\cdots +n_t$.  
\end{corollary 3.4}   

\begin{proof} Let $k$ be as in the Theorem 3.3.
By Theorem 1.1 (Compatibility), it suffices to show that
for each $\underline n \in \mathbb N^t$ and each $P \in
\ass(M/I_1^{n_1}I_2^{n_2}\cdots I_t^{n_t}M)$, there is a $P$-primary
component $Q$ of $I_1^{n_1}I_2^{n_2}\cdots I_t^{n_t}M \subset M$ such
that $P^{k|\underline n|}M \subseteq Q$. So we fix $\underline n$
and $P \in \ass(M/I_1^{n_1}I_2^{n_2}\cdots I_t^{n_t}M)$. Let 
\begin{eqnarray*}
(P^{k|\underline n|}M+I_1^{n_1}I_2^{n_2}\cdots I_t^{n_t}M)
&=& Q_1 \cap Q_2 \cap \cdots \cap Q_r \qquad \text {and} \\
(I_1^{n_1}I_2^{n_2}\cdots I_t^{n_t}M :_M  P^{\infty})
&=& Q_{r+1} \cap Q_{r+2} \cap \cdots \cap Q_s 
\end{eqnarray*}
be irredundant and minimal primary decompositions of the corresponding
submodules of $M$, where $Q_i$ is a $P_i$-primary submodule of $M$ for
each $i=1,2,\dots,s$. As $P \notin \ass(M/(I_1^{n_1}I_2^{n_2}\cdots
I_t^{n_t}M :_M  P^{\infty}))$, we may assume that $P_1=P$.
By Theorem 3.3,  
$(P^{k|\underline n|}M+I_1^{n_1}I_2^{n_2}\cdots I_t^{n_t}M) \cap
(I_1^{n_1}I_2^{n_2}\cdots I_t^{n_t}M :_M  P^{\infty})=
I_1^{n_1}I_2^{n_2}\cdots I_t^{n_t}M$. Hence 
\[
I_1^{n_1}I_2^{n_2}\cdots I_t^{n_t}M = Q_1 \cap Q_2 \cap \cdots \cap
Q_r \cap Q_{r+1} \cap Q_{r+2} \cap \cdots \cap Q_s.
\]
Although the above intersection may not necessarily be irredundant and
minimal, we know that $Q_1$ is a $P_1=P$-primary component of
$I_1^{n_1}I_2^{n_2}\cdots I_t^{n_t}M \subset M$ since $P \in
\ass(M/I_1^{n_1}I_2^{n_2}\cdots I_t^{n_t}M)$ and $Q_1$ is the only
$P$-primary submodule in the above intersection. Evidently
$P^{k|\underline n|}M \subseteq Q_1$.
\end{proof}

\addtocounter{theorem}{1}

Actually Theorem 3.3 can be stated in a more general situation:
The filtration $\{I_1^{n_1}I_2^{n_2}\cdots I_t^{n_t}M \,| \,(n_1, n_2,
\dots,n_t) \in \mathbb Z^t \}$ may be replaced by a `multi-graded'
filtration $\{M_{(n_1, n_2, \dots, n_t)} \,| \, (n_1, n_2, \dots,n_t)
\in \mathbb Z^t \}$ of $M$ such that  
\[
\mathcal{M}= \bigoplus_{(n_1, n_2, \dots , n_t) \in \mathbb Z^t}M_{(n_1, n_2,
\dots, n_t)}T_1^{n_1}T_2^{n_2}\cdots T_t^{n_t} 
\]
naturally forms a multi-graded Noetherian module over a multi-graded
sub-ring $\mathcal{R}$ in $R[T_1,T_2,\dotsc, T_t,
T_1^{-1},T_2^{-1},\dots,T_t^{-1}] $ with the usual grading such that
$T_1^{-1},T_2^{-1},\linebreak[0]\dotsc,T_t^{-1}$ are all contained in
$\mathcal{R}$ and 
the $(0,0,\dots,0)$ part of $\mathcal{R}$ is $R$. We
call such a filtration `Noetherian'. To simplify notation, we use
$\underline n$ to denote $(n_1, n_2, \dots, n_t)$ and use
$|\underline n |$ to denote $n_1+ n_2+ \dots + n_t$. And $\mathbb N^t :=
\{ (n_1, n_2, \dots, n_t) \,| \, n_i \ge 0, \ i=1,2,\dots, t\}$.

The next theorem and its corollary look apparently more general than
Theorem 3.3 and Corollary 3.4, although in essence they are the same.

\begin{theorem} Let $R$ be a Noetherian ring,  $M$ a finitely
generated $R$-module and $\{M_{(n_1 n_2 \dots n_t)} \,| \, (n_1, n_2,
\dots,n_t) \in \mathbb Z^t \}$ a Noetherian filtration of $M$. Then 
\begin{enumerate}
\item There
exists a $k \in \mathbb N$ such that for all $\underline m \in \mathbb Z^t$
, for all $\underline n \in \mathbb N^t$ and 
for all ideals $J \subset R$, $ J^{k|\underline n|}M_{\underline m}
\cap (M_{\underline m+\underline n} :_{M_{\underline m}}
J^{\infty}) \subseteq M_{\underline m+\underline n},$
i.e. $(J^{k|\underline n|}M_{\underline m}+M_{\underline
m+\underline n})  \cap
(M_{\underline m+\underline n} :_{M_{\underline m}} J^{\infty})=
M_{\underline m+\underline n}$;
\item The set 
$ \cup_{\underline m \in \mathbb Z^t,\underline n \in \mathbb N^t}\ass
(M_{\underline m}/M_{\underline m+\underline n}) $ is finite.
\end{enumerate}
\end{theorem} 
 
\begin{proof}The proof of (1) may be carried out in almost the same way
as in the proof of Theorem 3.3. But here we choose to use Theorem 3.3
and provide a sketch of the proof: Simply apply Theorem 3.3 to the
Noetherian $\mathcal{R}$ module $\mathcal{M}$ and ideals
$\mathcal{I}_i=T_i^{-1}\mathcal{R}$ and then restrict the results to each of
the homogeneous 
pieces. Theorem 3.3 gives results for all the ideals of $\mathcal{R}$, but
here we are only interested in the ideals $J\mathcal{R}$, the ideals
extended from ideals $J \subset R$.

To prove (2), we notice that the set
\[
\bigcup_{\underline n \in \mathbb N^t} \ass_{\mathcal{R}}
(\mathcal{M}/T_1^{-n_1} T_2^{-n_2} \cdots T_t^{-n_t}\mathcal{M})
\]
is finite. Then (2) follows by contracting to each of the homogeneous
pieces. 
\end{proof}


\begin{corollary}  Let $R$ be a Noetherian ring,  $M$ a finitely
generated $R$-module and $\{M_{(n_1 n_2 \dots n_t)} \,| \, (n_1, n_2,
\dots,n_t) \in \mathbb Z^t \}$ a Noetherian filtration of $M$. Then there
exists a $k \in \mathbb N$ such that for any $\underline m \in \mathbb Z^t$
, $\underline n \in \mathbb N^t$ and $P \in \ass(M_{\underline
  m}/M_{\underline m+\underline n})$, 
there exists a $Q \in \Lambda_P(M_{\underline m+\underline n} \subseteq
M_{\underline m} )$ such that
$P^{k|\underline n|}M_{\underline m} \subseteq Q$.
\end{corollary}

\begin{example} [Compare with \cite{Sh1}] Assume that $R$ is
 Nagata (e.g. $R$ is 
excellent) and $M$ is a finitely generated $R$-module and $I_1, I_2,
\dots, I_t$ ideals  of $R$. Then we have a multi-graded filtration
$\{\overline {I_1^{n_1}I_2^{n_2}\cdots I_t^{n_t}}M \,| \, \underline n \in
\mathbb Z^t \}$. In order to see if the filtration satisfies the Linear
Growth property, we may mod out the nil-radical and hence assume that
$R$ is reduced. Then it is straightforward to see that the associated
graded module is finite over $\mathcal{R}=R[I_1T_1,I_2T_2,\cdots, I_tT_t,
T_1^{-1},T_2^{-1},\dots,T_t^{-1}]$. Hence the filtration  satisfies
the Linear Growth property. Similarly we can show the Linear Growth
property of the filtration  $\{ \overline {I_1^{n_1}}\cdot \overline
{I_2^{n_2}} \cdots \overline {I_t^{n_t}}M \,| \, \underline n \in
\mathbb Z^t \}$ provided $R$ is reduced and Nagata.
\end{example}

In \cite{Sh1} R. Y. Sharp proved the Linear Growth property of the
filtration $\{\overline {I^n}\,|\, n \in \mathbb Z\}$ of Noetherian ring $R$
without the Nagata assumption. The argument there also works for the
filtration $ \{\overline {I_1^{n_1}I_2^{n_2}\cdots I_t^{n_t}} \,| \,
\underline n \in \mathbb Z^t \}$ of any Noetherian ring $R$. That is
because the set $\cup_{\underline n \in \mathbb Z^t}\ass (R/\overline
{I_1^{n_1}I_2^{n_2}\cdots I_t^{n_t}})$ is finite (cf. \cite{Ra}) and
hence we can \linebreak
localize and then complete. In fact, if we know in advance the set 
\linebreak 
$\cup_{\underline n \in \mathbb Z^t}\ass (M/\overline 
{I_1^{n_1}I_2^{n_2}\cdots I_t^{n_t}}M)$ is finite for a
finitely generated faithful $R$-module $M$, we can localize and then
complete and then contract the result of Example 3.7 for $\hat M$
back to $M$  to
deduce that the filtration $\{\overline {I_1^{n_1}I_2^{n_2}\cdots
I_t^{n_t}}M \,| \, \underline n \in \mathbb Z^t \}$ satisfies the Linear
Growth property. We need $M$ to be faithful so that the process of
contraction works.

\begin{example} Assume $R$ is Nagata and has characteristic
$p$, where $p$ is a prime number and $M$ is a finitely generated
$R$-module. Then for any  ideal $I$ in $R$,
tight closure of $I$, denoted by $I^*$, is defined \cite{HH}.
It is shown that $ \sqrt 0 \subseteq I^* \subseteq \overline I$ for any
ideal $I$ in $R$ \cite{HH}. By the
same argument as in Example 3.7 we can deduce that the filtration  
$\{ (I_1^{n_1}I_2^{n_2}   \cdots I_t^{n_t})^*M \,| \, \underline n \in
\mathbb Z^t \}$ is Noetherian and hence has the Linear Growth
property. If, furthermore, $R$ is 
reduced, then the filtration  $\{ I_1^{n_1*}I_2^{n_2*} \cdots
I_t^{n_t*}M \,| \, \underline n \in \mathbb Z^t \}$ satisfies the Linear
Growth property. 
\end{example}

In \cite{Ra} it is shown that $\ass(R/\overline {I^n})$ is
non-decreasing and eventually stabilizes for any ideal $I$ in a
Noetherian ring 
$R$. For any finitely generated $R$-module $M$, a result of \cite{Br}
says that $\ass (M/I^{n}M)$ also stabilizes for large $n$. If $R$ is
Nagata and of characteristic $p>0$, then it follows from Example
3.8 and Theorem 3.5 that the set $\cup_{\underline n \in \mathbb
Z^t}\ass (M/(I_1^{n_1}I_2^{n_2}\cdots I_t^{n_t})^*M)$ is finite. In
case of $t=1$, we would like to study the stability of $\ass
(M/(I^{n*}M))$. Since $\oplus_{n \in \mathbb Z}I^{n*}MT^n$ is finite over
$R[IT, T^{-1}]$ (see Example 3.8), we know the filtration $\{I^{n*}M
\,|\, n\in \mathbb N \}$ of $M$ is eventually stable, i.e. $I^{n+1*}M =
II^{n*}M$ for all large $n$. Hence the argument in \cite{Br} can be
applied to show that  $\ass (M/I^{n*}M)$ stabilizes for large $n$. 

\section*{Acknowledgment} 

I would like to thank Craig Huneke for the helpful conversations and
valuable inspiration, which made this paper possible.

\end{document}